\input amstex
\input amsppt.sty
\magnification=\magstep1
\hsize=30truecc
\vsize=22.2truecm
\baselineskip=16truept
\NoBlackBoxes
\TagsOnRight \pageno=1 \nologo
\def\Z{\Bbb Z}
\def\N{\Bbb N}

\def\l{\left}
\def\r{\right}
\def\bg{\bigg}
\def\({\bg(}
\def\[{\bg\lfloor}
\def\){\bg)}
\def\]{\bg\rfloor}
\def\t{\text}
\def\f{\frac}

\def\p{\ (\roman{mod}\ p)}

\def\bi{\binom}
\def\eq{\equiv}

\def\ls{\leqslant}
\def\gs{\geqslant}
\def\mo{\roman{mod}}

\def\da{\delta}

\def\Proof{\noindent{\it Proof}}

\def\Remark{\medskip\noindent{\it  Remark}}

\def\Ack{\medskip\noindent {\bf Acknowledgment}}
\hbox {Colloq. Math. 130(2013), no.\,1, 67--78.}
\bigskip
\topmatter
\title Arithmetic Theory of Harmonic Numbers (II)\endtitle
\author Zhi-Wei Sun and Li-Lu Zhao\endauthor
\leftheadtext{Zhi-Wei Sun and Li-Lu Zhao} \rightheadtext{Arithmetic theory of
harmonic numbers (II)}
\author Zhi-Wei Sun$^1$ and Li-Lu Zhao$^2$\endauthor
\leftheadtext{Zhi-Wei Sun and Li-Lu Zhao}
\affil $^1$Department of
Mathematics, Nanjing University\\
Nanjing 210093, People's Republic of China\\zwsun\@nju.edu.cn
\\{\tt http://math.nju.edu.cn/$\sim$zwsun}
\medskip
$^2$School of Mathematics
\\Hefei University of Technology
\\Hefei 230009, People's Republic of China
\\zhaolilu\@gmail.com
\endaffil
\abstract For $k=1,2,\ldots$ let $H_k$ denote the harmonic number
$\sum_{j=1}^k 1/j$. In this paper we establish some new congruences
involving harmonic numbers. For example, we show that for any prime
$p>3$ we have
$$\sum_{k=1}^{p-1}\f{H_k}{k2^k}\eq\f7{24}pB_{p-3}\pmod{p^2},\ \ \sum_{k=1}^{p-1}\f{H_{k,2}}{k2^k}\eq-\f 38B_{p-3}\pmod{p},$$
and
$$\sum_{k=1}^{p-1}\f{H_{k,2n}^2}{k^{2n}}\eq\f{\bi{6n+1}{2n-1}+n}{6n+1}pB_{p-1-6n}\ (\mo\ p^2)$$
for any positive integer $n<(p-1)/6$, where $B_0,B_1,B_2,\ldots$ are
Bernoulli numbers, and $H_{k,m}:=\sum_{j=1}^k 1/j^m$.
\endabstract
\thanks 2010 {\it Mathematics Subject Classification}.\,Primary 11A07, 11B68;
Secondary 05A19,\ 11B75.
\newline \indent The first author is supported by the National Natural Science
Foundation (grant 11171140) of China and the PAPD of Jiangsu Higher
Education Institutions.
\endthanks
\endtopmatter
\document

\heading{1. Introduction}\endheading

Recall that harmonic numbers are those
$$H_n:=\sum_{0<k\ls n}\f1k\ \ \ (n\in\N=\{0,1,2,\ldots\}),$$
where $H_0:=0$ since we consider the value of an empty sum as zero.
They play important roles in mathematics. In 1862 J. Wolstenholme
[W] showed the congruence $H_{p-1}\eq0\pmod{p^2}$ for any prime
$p>3$. Throughout this paper, for a prime $p$ and two rational
$p$-adic integers $A$ and $B$, we write $A\eq B\ (\mo\ p^n)$ (with
$n\in\N$) to mean that $A-B$ is divisible by $p^n$ in the ring of
$p$-adic integers.

In [Su] the first author investigated arithmetic properties of
harmonic numbers systematically. For example, he proved that for any
prime $p>5$ we have
$$\sum_{k=1}^{p-1}\f{H_k}{k2^k}\eq\sum_{k=1}^{p-1}\f{H_k^2}{k^2}\eq0\ (\mo\ p).$$

For $m\in\Z^+=\{1,2,3,\ldots\}$, harmonic numbers of order $m$ are
defined by
$$H_{n,m}:=\sum_{0<k\ls n}\f1{k^m}\ \ (n\in\N).$$

It is known that
$$\sum_{k=1}^\infty\f{H_k}{k2^k}=\f{\pi^2}{12}\quad \ \t{(S. W. Coffman [C], 1987)}$$
and
$$\sum_{k=1}^\infty\f{H_{k,2}}{k2^k}=\f 58\zeta(3)\quad \ \t{(B. Cloitre, 2004)}.$$
Both identities can be found in [SW].

Our first theorem is as follows.
\proclaim{Theorem 1.1} For any prime $p>3$, we have
$$\sum_{k=1}^{p-1}\f{H_k}{k2^k}\eq\f 7{24}pB_{p-3}\pmod{p^2}\tag1.1$$
and
$$\sum_{k=1}^{p-1}\f{H_{k,2}}{k2^k}\eq-\f 38B_{p-3}\ (\mo\ p),\tag1.2$$
where $B_0,B_1,B_2,\ldots$ are Bernoulli numbers.
\endproclaim
\Remark\ 1.1. (1.1) confirms the first part of [Su, Conjecture 1.1].
The second part of [Su, Conjecture
1.1] states that $\sum_{k=1}^{p-1}H_k^2/k^2\eq\f 45pB_{p-5}\pmod{p^2}$
for any prime $p>3$; this was confirmed by R. Me\v strovi\'c [M] quite recently.
\medskip

Our second theorem confirms the second conjecture of [Su].

\proclaim{Theorem 1.2 {\rm ([Su, Conjecture 1.2])}} Let $p$ be an odd prime and let $n$ be a
positive integer with  $p-1\nmid 6n$. Then
$$\sum_{k=1}^{p-1}\f{H_{k,2n}^2}{k^{2n}}\eq0\ (\mo\ p).\tag1.3$$
Furthermore, when $p>6n+1$ we have
$$\sum_{k=1}^{p-1}\f{H_{k,2n}^2}{k^{2n}}\eq\f{s(n)}{6n+1}pB_{p-1-6n}\ (\mo\ p^2),\tag1.4$$
where
$$s(n)=\bi{6n+1}{2n-1}+n.$$
\endproclaim
\Remark\ 1.2. We give here four initial values of the integer sequence $\{s(n)\}_{n\gs1}$:
$$s(1)=8,\ s(2)=288,\ s(3)=11631,\ s(4)=480704.$$
\medskip

We will show Theorems 1.1 and 1.2 in Sections 2 and 3 respectively.

\heading{2. Proof of Theorem 1.1}\endheading

\proclaim{Lemma 2.1} Let $p>3$ be a prime. Then
$$\sum_{k=1}^{p-1}\f{(-1)^k}{k^2}\eq\f p2B_{p-3}\pmod{p^2},\ \sum_{k=1}^{p-1}\f{(-1)^k}{k^3}\eq -\f{B_{p-3}}2\pmod{p},\tag2.1$$
and
$$\sum_{k=1}^{p-1}\f{H_k}k\eq\f p3B_{p-3}\pmod{p^2}\ \ \t{and}\ \ \sum_{k=1}^{p-1}\f{(-1)^k}{k^2}H_k\eq -\f{B_{p-3}}4\pmod{p}.\tag2.2$$
\endproclaim
\Proof. It is known that (cf. [S, Corollaries 5.1 and 5.2])
$$\sum_{k=1}^{p-1}\f1{k^2}\eq\f23pB_{p-3}\pmod{p^2},\ \ \sum_{k=1}^{p-1}\f1{k^3}\eq \f 34pB_{p-4}\eq-p\da_{p,5}\pmod{p^2},$$
and
$$\sum_{k=1}^{(p-1)/2}\f1{k^2}\eq\f 73pB_{p-3}\pmod{p^2}\ \ \t{and}\ \ \sum_{k=1}^{(p-1)/2}\f1{k^3}\eq-2B_{p-3}\pmod{p}.$$
Thus
$$\align \sum_{k=1}^{p-1}\f{(-1)^k}{k^2}=&\sum_{k=1}^{p-1}\f{1+(-1)^k}{k^2}-\sum_{k=1}^{p-1}\f1{k^2}=\f12H_{(p-1)/2,2}-H_{p-1,2}
\\\eq&\f 76pB_{p-3}-\f 23pB_{p-3}=\f p2B_{p-3}\pmod{p^2}
\endalign$$
and
$$\align\sum_{k=1}^{p-1}\f{(-1)^k}{k^3}=&\sum_{k=1}^{p-1}\f{1+(-1)^k}{k^3}-\sum_{k=1}^{p-1}\f1{k^3}
\\=&\f14H_{(p-1)/2,3}-H_{p-1,3}\eq\f{-2B_{p-3}}4\pmod{p}.
\endalign$$
Therefore (2.1) holds.

By the proof of [S, Theorem 6.1],
$$\sum_{1\ls j<k\ls p-1}\f1{jk}\eq-\f p3B_{p-3}\pmod{p^2}.$$
So we have
$$\sum_{k=1}^{p-1}\f{H_k}k=\sum_{k=1}^{p-1}\f1{k^2}+\sum_{1\ls j<k\ls p-1}\f1{jk}
\eq\f23pB_{p-3}-\f p3B_{p-3}=\f p3B_{p-3}\pmod{p^2}.$$
This proves the first congruence in (2.2).

Now we prove the second congruence in (2.2). Since
$$H_{p-k}=H_{p-1}-\sum_{j=1}^{k-1}\f1{p-j}\eq H_{k-1}=H_k-\f1k\pmod p$$
for all $k=1,\ldots,p-1$, we have
$$\sum_{k=1}^{p-1}\f{(-1)^k}{k^2}H_k=\sum_{k=1}^{p-1}\f{(-1)^{p-k}}{(p-k)^2}H_{p-k}
\eq-\sum_{k=1}^{p-1}\f{(-1)^k}{k^2}\l(H_k-\f1k\r)\pmod p$$
and hence
$$\sum_{k=1}^{p-1}\f{(-1)^k}{k^2}H_k\eq\f12\sum_{k=1}^{p-1}\f{(-1)^k}{k^3}\eq-\f{B_{p-3}}4\pmod p.$$

The proof of Lemma 2.1 is now complete. \qed

\proclaim{Lemma 2.2} {\rm (i)} For any positive integers $k$ and $m$ we have
$$\sum_{n=1}^m\bi{n-1}{k-1}=\bi mk.\tag2.3$$

{\rm (ii)} For each $n=1,2,3,\ldots$ we have
$$\sum_{k=1}^n\bi nk\f{(-1)^{k-1}}kH_k=H_{n,2}.\tag2.4$$
\endproclaim
\Proof. (2.3) is well known (cf. [G, (1.5)]) and it can be easily proved by induction on $m$.

(2.4) is also known (cf. [H]). Here we prove it by induction.
Clearly (2.4) holds for $n=1$. Assume that (2.4) holds for a fixed
positive integer $n$. Then
$$\align\sum_{k=1}^{n+1}\bi{n+1}k\f{(-1)^{k-1}}kH_k=&\sum_{k=1}^n\bi nk\f{(-1)^{k-1}}kH_k+\sum_{k=1}^{n+1}\bi n{k-1}\f{(-1)^{k-1}}kH_k
\\=&H_{n,2}+\f1{n+1}\sum_{k=0}^{n+1}\bi{n+1}k(-1)^{k-1}H_k.
\endalign$$
Note that
$$\align&\sum_{k=0}^{n+1}\bi{n+1}k(-1)^{k-1}H_k
\\=&\sum_{k=0}^{n}\bi{n}k(-1)^{k-1}H_k+\sum_{k=1}^{n+1}\bi{n}{k-1}(-1)^{k-1}\l(H_{k-1}+\f1k\r)
\\=&\sum_{k=1}^{n+1}\bi{n}{k-1}\f{(-1)^{k-1}}k=-\f1{n+1}\sum_{k=1}^{n+1}\bi {n+1}k(-1)^k=\f1{n+1}.
\endalign$$
So
$$\sum_{k=1}^{n+1}\bi{n+1}k\f{(-1)^{k-1}}kH_k=H_{n,2}+\f1{n+1}\cdot\f1{n+1}=H_{n+1,2}$$
as desired. \qed

\proclaim{Lemma 2.3} Let $p>3$ be a prime. Then
$$\sum_{1\ls j\ls k\ls p-1}\f{2^j(j+k)}{j^2k^2}\eq\sum_{k=1}^{p-1}\f{(-1)^k}{k^3}\pmod{p}.\tag2.5$$
\endproclaim
\Proof. Observe that
$$\align &\sum_{1\ls i\ls j\ls k\ls p-1}\f{2^i}{ijk}-\sum_{1\ls i<j<k\ls p-1}\f{2^i}{ijk}
\\=&\sum_{1\ls j\ls k\ls p-1}\f{2^j}{j^2k}+\sum_{1\ls i\ls j\ls p-1}\f{2^i}{ij^2}-\sum_{k=1}^{p-1}\f{2^k}{k^3}
\\=&\sum_{1\ls j\ls k\ls p-1}\l(\f{2^j}{j^2k}+\f{2^j}{jk^2}\r)-\sum_{k=1}^{p-1}\f{2^k}{k^3}.
\endalign$$
Similarly,
$$\align &2\sum_{1\ls i\ls j\ls k\ls p-1}\f{(-1)^i}{ijk}-2\sum_{1\ls i<j<k\ls p-1}\f{(-1)^i}{ijk}-2\sum_{k=1}^{p-1}\f{(-1)^k}{k^3}
\\=&2\sum_{1\ls j<k\ls p-1}\l(\f{(-1)^j}{j^2k}+\f{(-1)^j}{jk^2}\r)
\\\eq&\sum_{1\ls j<k\ls p-1}\l(\f{(-1)^j}{j^2k}+\f{(-1)^j}{jk^2}+\f{(-1)^{p-j}}{(p-j)^2(p-k)}+\f{(-1)^{p-j}}{(p-j)(p-k)^2}\r)
\\=&\sum_{1\ls j<k\ls p-1}\f{(-1)^j}{j^2k}+\sum_{1\ls k<j\ls p-1}\f{(-1)^j}{j^2k}
\\&+\sum_{1\ls j<k\ls p-1}\f{(-1)^j}{jk^2}+\sum_{1\ls k<j\ls p-1}\f{(-1)^j}{jk^2}
\\=&H_{p-1}\sum_{j=1}^{p-1}\f{(-1)^j}{j^2}+H_{p-1,2}\sum_{j=1}^{p-1}\f{(-1)^j}j-2\sum_{j=1}^{p-1}\f{(-1)^j}{j^3}\pmod{p}.
\endalign$$
Thus, with the help of  $H_{p-1}\eq H_{p-1,2}\eq0\pmod p$, we have
$$\sum_{1\ls i\ls j\ls k\ls p-1}\f{(-1)^i}{ijk}\eq\sum_{1\ls i<j<k\ls p-1}\f{(-1)^i}{ijk}\pmod{p}.$$
By [ZS, Theorem 1.2],
$$\sum_{1\ls i<j<k\ls p-1}\f{(1-x)^i}{ijk}\eq\sum_{1\ls i<j<k\ls p-1}\f {x^i}{ijk}\pmod{p}.$$
So, in view of the above, we have
$$\align\sum_{1\ls i\ls j\ls k\ls p-1}\f{(-1)^i}{ijk}\eq&\sum_{1\ls i<j<k\ls p-1}\f {2^i}{ijk}
\\\eq&\sum_{1\ls i\ls j\ls k\ls p-1}\f{2^i}{ijk}+\sum_{k=1}^{p-1}\f{2^k}{k^3}
-\sum_{1\ls j\ls k\ls p-1}\f{2^j(j+k)}{j^2k^2}\pmod{p}.
\endalign$$
It remains to show that
$$\sum_{1\ls i\ls j\ls k\ls p-1}\f{2^i-(-1)^i}{ijk}\eq\sum_{k=1}^{p-1}\f{(-1)^k-2^k}{k^3}\pmod p.\tag2.6$$
With the help of Lemma 2.2, we have
$$\align\sum_{1\ls i\ls j\ls k\ls p-1}\f{2^i-(-1)^i}{ijk}=&\sum_{1\ls i\ls j\ls k\ls p-1}\f1{ijk}\sum_{r=0}^{i}(1-(-2)^r)\bi ir
\\=&\sum_{r=1}^{p-1}\f{1-(-2)^r}r\sum_{1\ls  j\ls k\ls p-1}\f1{jk}\sum_{i=1}^j\bi{i-1}{r-1}
\\=&\sum_{r=1}^{p-1}\f{1-(-2)^r}r\sum_{1\ls  j\ls k\ls p-1}\f1{jk}\bi jr
\\=&\sum_{r=1}^{p-1}\f{1-(-2)^r}{r^2}\sum_{k=1}^{p-1}\f1k\sum_{j=1}^k\bi{j-1}{r-1}
\\=&\sum_{r=1}^{p-1}\f{1-(-2)^r}{r^2}\sum_{k=1}^{p-1}\f1k\bi kr=\sum_{r=1}^{p-1}\f{1-(-2)^r}{r^3}\sum_{k=1}^{p-1}\bi{k-1}{r-1}
\\=&\sum_{r=1}^{p-1}\f{1-(-2)^r}{r^3}\bi {p-1}r\eq\sum_{r=1}^{p-1}\f{(-1)^r-2^r}{r^3}\pmod p.
\endalign$$
This proves the desired (2.6). \qed

\medskip
\noindent{\it Proof of Theorem 1.1}. We prove (1.2) first. In view of (2.4), we have
$$\align \sum_{n=1}^{p-1}\f{H_{n,2}}{n2^n}=&\sum_{n=1}^{p-1}\f1{n2^n}\sum_{k=1}^n\bi nk\f{(-1)^{k-1}}kH_k
\\=&\sum_{k=1}^{p-1}\f{(-1)^{k-1}}kH_k\sum_{n=k}^{p-1}\f1{n2^n}\bi nk=\sum_{k=1}^{p-1}\f{(-1)^{k-1}}{k^22^k}H_k\sum_{n=k}^{p-1}\bi{n-1}{k-1}\f1{2^{n-k}}
\\=&\sum_{k=1}^{p-1}\f{(-1)^{k-1}}{k^22^k}H_k\sum_{j=0}^{p-1-k}\bi{k+j-1}j\f1{2^j}
\\=&\sum_{k=1}^{p-1}\f{(-1)^{k-1}}{k^22^k}H_k\sum_{j=0}^{p-1-k}\bi{-k}j\f1{(-2)^j}
\endalign$$
and hence
$$\align \sum_{n=1}^{p-1}\f{H_{n,2}}{n2^n}\eq&\sum_{k=1}^{p-1}\f{(-1)^{k-1}}{k^22^k}H_k\sum_{j=0}^{p-1-k}\bi{p-k}j\f1{(-2)^j}
\\=&\sum_{k=1}^{p-1}\f{(-1)^{k-1}}{k^22^k}H_k\f{1+(-1)^k}{2^{p-k}}
\\\eq&-\f12\sum_{k=1}^{p-1}\f{H_k}{k^2}\l(1+(-1)^k\r)\pmod{p}.
\endalign$$
Note that
$$\sum_{k=1}^{p-1}\f{H_k}{k^2}\eq B_{p-3}\ (\mo\ p)\ \ \t{and}\ \
\sum_{k=1}^{p-1}\f{(-1)^k}{k^2}H_k\eq-\f{B_{p-3}}4\ (\mo\ p)$$
by [ST, (5.4)] and (2.2) respectively. So we get
$$\sum_{n=1}^{p-1}\f{H_{n,2}}{n2^n}\eq-\f12\l(B_{p-3}-\f{B_{p-3}}4\r)=-\f 38B_{p-3}\pmod{p}.$$

Now we show (1.1).
Observe that
$$\align\sum_{k=1}^{p-1}\f{H_k}{k2^k}=&\sum_{1\ls j\ls k\ls p-1}\f1{jk2^k}=\sum_{1\ls j\ls k\ls p-1}\f1{(p-k)(p-j)2^{p-j}}
\\=&\sum_{1\ls j\ls k\ls p-1}\f{2^{j-p}(p+j)(p+k)}{(p^2-j^2)(p^2-k^2)}
\\\eq&\sum_{1\ls j\ls k\ls p-1}\f{2^{j-p}(jk+p(j+k))}{j^2k^2}
\\\eq&2^{-p}\sum_{1\ls j\ls k\ls p-1}\f{2^j}{jk}
+\f p2\sum_{1\ls j\ls k\ls p-1}\f{2^j(j+k)}{j^2k^2}\pmod{p^2}.
\endalign$$
In view of Lemmas 2.2 and 2.1,
$$\align \sum_{1\ls j\ls k\ls p-1}\f{2^j-1}{jk}=&\sum_{1\ls j\ls k\ls p-1}\f1{jk}\sum_{i=1}^j\bi ji
=\sum_{i=1}^{p-1}\f1i\sum_{k=1}^{p-1}\f1k\sum_{j=1}^k\bi{j-1}{i-1}
\\=&\sum_{i=1}^{p-1}\f1i\sum_{k=1}^{p-1}\f1k\bi ki=\sum_{i=1}^{p-1}\f1{i^2}\sum_{k=1}^{p-1}\bi{k-1}{i-1}
\\=&\sum_{i=1}^{p-1}\f1{i^2}\bi{p-1}i
=\sum_{i=1}^{p-1}\f{(-1)^i}{i^2}\prod_{r=1}^i\l(1-\f pr\r)
\\\eq&\sum_{i=1}^{p-1}\f{(-1)^i(1-pH_i)}{i^2}\eq\f p2B_{p-3}-p\l(-\f{B_{p-3}}4\r)\pmod{p^2}.
\endalign$$
Note that
$$\sum_{1\ls j\ls k\ls p-1}\f1{jk}=\sum_{k=1}^{p-1}\f{H_k}k\eq\f p3B_{p-3}\pmod{p^2}$$
by (2.2). Combining the above with (2.5), we finally obtain that
$$\align\sum_{k=1}^{p-1}\f{H_k}{k2^k}\eq& 2^{-p}\l(\f 34pB_{p-3}+\f p3B_{p-3}\r)+\f p2\sum_{k=1}^{p-1}\f{(-1)^k}{k^3}
\\\eq&\f {13}{24}pB_{p-3}+\f p2\l(-\f{B_{p-3}}2\r)=\f 7{24}pB_{p-3}\pmod{p^2}\ \ (\t{by (2.1)}).
\endalign$$
This concludes the proof. \qed

\heading{3. Proof of Theorem 1.2}\endheading

\proclaim{Lemma 3.1} Let $p>3$ be a prime and let $m$ be a positive integer with $p-1\nmid 3m$. Then
$$\sum_{1\ls j<k\ls p-1}\l(\f1{j^mk^{2m}}+\f1{j^{2m}k^m}\r)\eq 0\pmod{p}.\tag3.1$$
Moreover, if $p>3m+1$, then
$$\sum_{1\ls j<k\ls p-1}\l(\f1{j^mk^{2m}}+\f1{j^{2m}k^m}\r)\eq-p\f{3m}{3m+1}B_{p-1-3m}\pmod{p^2}.\tag3.2$$
\endproclaim
\Proof. It is well-known that
$$\sum_{k=1}^{p-1}\f1{k^n}\eq0\pmod{p}\quad\t{for any integer}\ n\not\eq0\pmod{p-1}.$$
Also,
$$\sum_{k=1}^{p-1}\f1{k^n}\eq\f {pn}{n+1}B_{p-1-n}\pmod{p^2}\quad\t{for}\ n=1,\ldots,p-2$$
(see, e.g., [S, Corollary 5.1]). Thus
$$\sum_{1\ls j<k\ls p-1}\l(\f1{j^mk^{2m}}+\f1{j^{2m}k^m}\r)
=\sum_{j=1}^{p-1}\f1{j^m}\sum_{k=1}^{p-1}\f1{k^{2m}}-\sum_{k=1}^{p-1}\f1{k^{3m}}\eq0\pmod p.$$
Moreover, we have (3.2) if $p>3m+1$. \qed

\proclaim{Lemma 3.2} Let $p>3$ be a prime and let $m$ be a positive even integer.
Then
$$\sum_{1\ls j<k\ls p-1}\l(\f1{j^mk^{2m}}-\f1{j^{2m}k^m}\r)\eq 0\pmod{p}.\tag3.3$$
Moreover, if $p>3m+1$ then
$$\sum_{1\ls j<k\ls p-1}\l(\f1{j^mk^{2m}}-\f1{j^{2m}k^m}\r)\eq \f {pm\bi{3m}mB_{p-1-3m}}{(m+1)(2m+1)}\ \ (\mo\ p^2).\tag3.4$$
\endproclaim
\Proof. As $m$ is even, we have
$$\align\sum_{1\ls j<k\ls p-1}\f1{j^mk^{2m}}=&\sum_{1\ls j<k\ls p-1}\f1{(p-k)^m(p-j)^{2m}}
\\\eq&\sum_{1\ls j<k\ls p-1}\f1{j^{2m}k^m}\pmod{p}.
\endalign$$

Now suppose that $p>3m+1$. Then
$$\align&\sum_{1\ls j<k\ls p-1}\f1{j^mk^{2m}}=\sum_{1\ls j<k\ls p-1}\f{(p+k)^m(p+j)^{2m}}{(p^2-k^2)^m(p^2-j^2)^{2m}}
\\\eq&\sum_{1\ls j<k\ls p-1}\f{(k^m+pmk^{m-1})(j^{2m}+p2mj^{2m-1})}{j^{4m}k^{2m}}
\\\eq&\sum_{1\ls j<k\ls p-1}\f1{j^{2m}k^m}+pm\sum_{1\ls j<k\ls p-1}\l(\f1{j^{2m}k^{m+1}}+\f2{j^{2m+1}k^m}\r)\pmod{p^2}.
\endalign$$
So, (3.4) is reduced to
$$\sum_{1\ls j<k\ls p-1}\l(\f1{j^{2m}k^{m+1}}+\f2{j^{2m+1}k^m}\r)\eq\f{\bi{3m}mB_{p-1-3m}}{(m+1)(2m+1)}\pmod{p}.\tag3.5$$
Recall that for any integer $n$ we have
$$\sum_{k=1}^{p-1}k^n\eq\cases p-1\pmod{p}&\t{if}\ p-1\mid n,\\0\pmod p&\t{if}\ p-1\nmid n.\endcases$$
(See, e.g., [IR, p.235].) Also,
$$\sum_{j=0}^{k-1}j^n=\f1{n+1}\sum_{j=0}^n\bi{n+1}jB_jk^{n+1-j}$$
for any $k=1,2,3\ldots$ and $n=0,1,2,\ldots$.
(See, e.g., [IR, p.\,230].)
Therefore
$$\align &\sum_{1\ls j<k\ls p-1}\f1{j^{2m}k^{m+1}}
\\\eq&\sum_{k=1}^{p-1}\f1{k^{m+1}}\sum_{j=0}^{k-1}j^{p-1-2m}
=\sum_{k=1}^{p-1}\f1{k^{m+1}(p-2m)}\sum_{j=0}^{p-1-2m}\bi{p-2m}jB_jk^{p-2m-j}
\\\eq&-\f1{2m}\sum_{j=0}^{p-1-2m}\bi{p-2m}jB_j\sum_{k=1}^{p-1}k^{p-1-3m-j}
\\\eq&\f1{2m}\sum^{p-1-2m}\Sb j=0\\p-1\mid j+3m\endSb\bi{p-2m}jB_j=\f1{2m}\bi{p-2m}{m+1}B_{p-1-3m}
\\\eq&\f1{2m}\bi{-2m}{m+1}B_{p-1-3m}=\f{(-1)^{m+1}}{2m}\bi{3m}{m+1}B_{p-1-3m}\pmod p.
\endalign$$
Similarly,
$$\align &\sum_{1\ls j<k\ls p-1}\f1{j^{2m+1}k^{m}}
\\\eq&\sum_{k=1}^{p-1}\f1{k^{m}}\sum_{j=0}^{k-1}j^{p-2-2m}
=\sum_{k=1}^{p-1}\f1{k^{m}(p-1-2m)}\sum_{j=0}^{p-2-2m}\bi{p-1-2m}jB_jk^{p-1-2m-j}
\\\eq&-\f1{2m+1}\sum_{j=0}^{p-2-2m}\bi{p-1-2m}jB_j\sum_{k=1}^{p-1}k^{p-1-3m-j}
\\\eq&\f1{2m+1}\sum^{p-2-2m}\Sb j=0\\p-1\mid j+3m\endSb\bi{p-1-2m}jB_j
=\f1{2m+1}\bi{p-1-2m}{m}B_{p-1-3m}
\\\eq&\f1{2m+1}\bi{-1-2m}{m}B_{p-1-3m}
=\f{(-1)^{m}}{2m+1}\bi{3m}{m}B_{p-1-3m}\pmod p.
\endalign$$
Therefore
$$\align&\sum_{1\ls j<k\ls p-1}\l(\f1{j^{2m}k^{m+1}}+\f2{j^{2m+1}k^m}\r)
\\\eq&\l(\f{(-1)^{m+1}}{2m}\bi{3m}{m+1}+2\f{(-1)^m}{2m+1}\bi{3m}m\r)B_{p-1-3m}
\\=&\f{(-1)^m}{(m+1)(2m+1)}\bi{3m}mB_{p-1-3m}\pmod p.
\endalign$$
So (3.5) holds as $m$ is even. \qed

\medskip
\noindent{\it Proof of Theorem 1.2}. Let $m=2n$. Clearly
$$\align\sum_{k=1}^{p-1}\f{H_{k,m}^2}{k^m}=&\sum_{k=1}^{p-1}\f1{k^m}\(\sum_{j=1}^k\f1{j^m}\)^2
\\=&\sum_{k=1}^{p-1}\f1{k^m}\(\sum_{j=1}^k\f1{j^{2m}}+2\sum_{1\ls i<j\ls k}\f1{i^mj^{m}}\)
\\=&H_{p-1,3m}+\sum_{1\ls j<k\ls p-1}\f1{j^{2m}k^m}+2\sum_{1\ls i<j\ls p-1}\f1{i^mj^{2m}}
\\&+2\sum_{1\ls i<j<k\ls p-1}\f1{i^mj^mk^m}
\endalign$$
and
$$\align H_{p-1,m}^3=&\sum_{i=1}^{p-1}\f1{i^m}\(\sum_{k=1}^{p-1}\f1{k^{2m}}+2\sum_{1\ls j<k\ls p-1}\f1{j^mk^m}\)
\\=&H_{p-1,3m}+3\sum_{1\ls j<k\ls p-1}\l(\f1{j^{2m}k^m}+\f1{j^mk^{2m}}\r)+6\sum_{1\ls i<j<k\ls p-1}\f1{i^mj^mk^m}.
\endalign$$
As $H_{p-1,m}\eq0\pmod p$, from the above we obtain
$$\align \sum_{k=1}^{p-1}\f{H_{k,m}^2}{k^m}\eq&H_{p-1,3m}+\sum_{1\ls j<k\ls p-1}\l(\f1{j^{2m}k^m}+\f2{j^mk^{2m}}\r)
\\&-\f{H_{p-1,3m}}3-\sum_{1\ls j<k\ls p-1}\l(\f1{j^{2m}k^m}+\f1{j^mk^{2m}}\r)
\\=&\f23H_{p-1,3m}+\sum_{1\ls j<k\ls p-1}\f1{j^mk^{2m}}\pmod{p^2}.
\endalign$$
Thus, by (3.1), (3.3) and the congruence $H_{p-1,3m}\eq0\pmod p$, we immediately get (1.3).

Below we assume that $p>3m+1$. Adding (3.2) and (3.4) we obtain
$$\align2\sum_{1\ls j<k\ls p-1}\f1{j^mk^{2m}}\eq &pmB_{p-1-3m}\l(-\f 3{3m+1}+\f{\bi{3m}m}{(m+1)(2m+1)}\r)
\\=&\f{pm}{3m+1}\l(\f{\bi{3m+1}m}{m+1}-3\r)B_{p-1-3m}\pmod{p^2}.
\endalign$$
Note also that
$$H_{p-1-3m}\eq p\f{3m}{3m+1}B_{p-1-3m}\pmod{p^2}.$$
Therefore
$$\align \sum_{k=1}^{p-1}\f{H_{k,m}^2}{k^m}\eq&\f 23\cdot p\f{3m}{3m+1}B_{p-1-3m}+\l(\f{\bi{3m+1}m}{m+1}-3\r)\f{pm/2}{3m+1}B_{p-1-3m}
\\=&\l(\f{\bi{3m+1}m}{m+1}+1\r)\f{pm/2}{3m+1}B_{p-1-3m}
\\=&\l(\bi{3m+1}{m-1}+\f m2\r)\f{pB_{p-1-3m}}{3m+1}\pmod{p^2}.
\endalign$$
This proves (1.4).

So far we have completed the proof of Theorem 1.2. \qed

\Ack. The authors wish to thank the referee for helpful comments.

\enddocument